\documentclass[a4paper,10pt]{amsart}
\title[Lie Symmetry and Numerical Analysis of KdV Equation]%
{Lie Symmetry Classification and Numerical Analysis of KdV Equation with Power-law Nonlinearity}
\author{Rohollah Bakhshandeh-Chamazkoti$^1$}
\author{Mohsen Alipour$^2$}

\thanks{$^{1,2}$Department of Mathematics, Faculty of Basic science,
Babol  University of Technology, Babol, Iran,\\
Email: $^1$r\_bakhshandeh@nit.ac.ir, $^2$m.alipor@nit.ac.ir}
\date{\today}
\usepackage[usenames,dvipsnames]{color}
\usepackage[latin1]{inputenc}
\usepackage{etex}
\usepackage{tikz}
\usetikzlibrary{matrix,arrows}
\usepackage{longtable}
\usepackage{ragged2e}
\usepackage{varioref}
\usepackage{amsfonts}
\usepackage{amsmath}
\usepackage{amssymb}
\usepackage{amsthm}
\usepackage{booktabs}
\usepackage[all]{xy}
\usepackage{fontenc}
\usepackage{ifthen}
\usepackage{braket}
\usepackage[pdftex]{hyperref}
\usepackage{amsmath}
\usepackage{amssymb}
\usepackage[mathscr]{eucal}
\usepackage[usenames]{color}

\newtheorem{thm}{Theorem}

\theoremstyle{remark}

\numberwithin{equation}{section}
\newcommand{\p}{\partial}

\newcommand{\X}{{\bf X}}

\newcommand{\be}{\begin{eqnarray}}
\newcommand{\ee}{\end{eqnarray}}
\newcommand{\bee}{\begin{eqnarray*}}
\newcommand{\eee}{\end{eqnarray*}}
\begin{document}
\begin{abstract}
In this paper, a complete  Lie
symmetry analysis of the damped wave equation with time-dependent coefficients
 is investigated. Then the invariant
solutions and the exact solutions generated from the symmetries are presented. Moreover,
a Lie algebraic classifications and the optimal system are discussed. Finally, using Chebyshev pseudo-spectral method (CPSM),
a numerical analysis  to solve the
invariant solutions corresponded the Lie symmetries of main equation is presented. This method applies the Chebyshev-Gauss-Lobatto points
as collocation points.
\end{abstract}
\maketitle
\section{Introduction}
The symmetry group analysis plays an critical role in the analysis of differential equations. The first paper on group classification methods is \cite{[Lie]}, where Lie proves that
a linear two-dimensional second-order partial differential equation may admit at most a three-parameter invariance group.
He computed the maximal invariance group of the one--dimensional heat conductivity equation and
applied the symmetries to construct invariant solutions.
The symmetry reduction is an interesting method for
solving nonlinear partial differential equations, \cite{[5],[6],[7]}. There have been some new generalizations of the classical
Lie group analysis for symmetry reductions. For instance, L. V. Ovsiannikov \cite{[4]} is one of the mathematicians which  extended the method of partially invariant solutions. His
works is based on the concept of an equivalence group, which is a Lie transformation group acting in the extended space
that is called jet space, and preserving the class of given  partial differential equations.

The nonlinear evolution equations are especially generated as a
extended  kinds of the well known equations like the Korteweg-de
Vries (KdV) equations and Kadomtsev-Petviasvilli equations, \cite{[1]}.

The KdV equation, with power law nonlinearity and
linear damping with dispersion has the following form
\be\label{eq:maineq}
u_t+[a(t)u^n-l(t)]u_x+b(t)u_{xxx}-c(t)u=0,
\ee
 where $a, b, c$ and $l$ are arbitrary smooth
functions with respect to $t$.
In \cite{[1]}, an exact solitary wave
solution of the KdV equation with
power law nonlinearity with time-dependent coefficients
of the nonlinear as well as the dispersion terms is obtained.

In \cite{[Khalique]}, the authors have investigated a exact solutions and Lie symmetries of the mKdV equation with time-dependent coefficients
of (\ref{eq:maineq}), in special cases of $\displaystyle{a(t)={1\over t}}$, $\displaystyle{b(t)={K\over t^2}}$, and
$\displaystyle{a(t)=K_0}$, $\displaystyle{b(t)=K_1\exp(-2K_0t)}$ where $K, K_0$ and  $K$ are constants. Gungor and et. al. have investigated a Lie symmetry classification
of this KdV equation, \cite{Gu}.

This paper is devoted  to calculate the symmetries
of  (\ref{eq:maineq}) equation. In the (\ref{eq:maineq}) equation,
$u_t$ shows the evolution term,
 $a(t)u^nu_x$ is the power law nonlinearity,
while $n$ is the index of power law and $b(t)u_{xxx}$
is the dispersion term. Moreover, $c(t)u$ is the linear damping while  $l(t)u_x$ is
the first order dispersion term, \cite{[1]}. In case of
$n = 1$, the  (\ref{eq:maineq}) equation is the KdV equation while for $n = 2$, we have the modified version of  KdV equation. As an example, a special  form of (\ref{eq:maineq}) with $a(t) = 6, n = 1,
b(t) = 1$ with $c(t) = 0$ and $l(t)=-1/2t$ in investigated. Lie symmetry analysis, invariant solution and
optimal system of this case is obtained.
%
\section{Lie Symmetry Methods}
\subsection{Preliminaries}
Consider a partial differential equation with  $p$ independent
variables and $q$ dependent variables with the one-parameter Lie group
of transformations
\be
x_i\longmapsto x_i+\epsilon\xi^i(x,u)+O(\epsilon^2),\qquad
u_{\alpha}\longmapsto
u_{\alpha}+\epsilon\varphi^{\alpha}(x,u)+O(\epsilon^2),\label{eq:3}
\ee
where $i=1,\ldots,p$ and $\alpha=1,\ldots,q$.  The general vector field
\be
\X=\sum_{i=1}^p\xi^i(x,u)\frac{\partial}{{\partial}x_i}+
\sum_{\alpha=1}^q\varphi^{\alpha}(x,u)\frac{\partial}{{\partial}u_{\alpha}},\label{eq:4}
\ee
on the $(x,u)$ space is given. So the characteristic of the vector field $\X$
is equal to
\be
Q^{\alpha}(x,u^{(1)})=\varphi^{\alpha}(x,u)-
\sum_{i=1}^p\xi^i(x,u)\frac{{\partial}u^{\alpha}}{{\partial}x_i},\;\;\;\;\;\;
\alpha=1,\ldots,q.\label{eq:5}
\ee
\begin{thm}\label{thm:1}
\cite{[3]} Let $\X$ be a vector field given by (\ref{eq:4}), and let
$Q=(Q^1,\ldots,Q^q)$ be its characteristic, as in (\ref{eq:5}).
The $n$-th prolongation of $\X$ is given explicitly by
\be
{\rm Pr}^{(n)}\X=\sum_{i=1}^p\xi^i(x,u)\frac{\partial}{{\partial}x^i}+\sum_{\alpha=1}^q\sum_J
\varphi_J^{\alpha}(x,u^{(n)})\frac{\partial}{{\partial}u_J^{\alpha}},\label{eq:6}
\ee
with coefficients
\be
\varphi_{J,i}^{\alpha}=D_i\varphi_J^{\alpha}-\sum_{j=1}^pD_i\xi^ju_{J,j}^{\alpha}~.\label{eq:7}
\ee
Here, $J=(j_1,\ldots,j_k)$, with $1\leq k\leq p$
is a multi-indices, and $D_i$ is total derivative and
subscripts of $u$ are derivative with respect to the respective
coordinates.
\end{thm}
\begin{thm}\label{thm:2}
 \cite{[3]} A connected group of transformations $G$ is a symmetry group of a
differential equation $\Delta=0$ if and only if the classical
infinitesimal symmetry condition
\begin{eqnarray}\label{eq:inv}
{\rm Pr}^{(n)}\X(\Delta)=0 \qquad {\rm whenever} \qquad \Delta=0,
\end{eqnarray}
holds for every infinitesimal generator $X\in\mathfrak{g}$ of $G$.
\end{thm}
%
%
\subsection{Governing Equation}
In order to find Lie point symmetries of the partial differential equation (\ref{eq:maineq}), we consider
one-parameter Lie group of transformations
\be\nonumber
\overline{x}&=& \xi(x, t, u, \epsilon),\\
\overline{t}&=& \tau(x, t, u, \epsilon),\\\label{eq:transfoerm}
\overline{u}&=& \varphi(x, t, u, \epsilon),\nonumber
\ee
under which (\ref{eq:maineq}) must be invariant. The group action is infinitesimally given by
\be\nonumber
\overline{x}&=& x+\epsilon\xi(x, t, u)+O(\epsilon^2),\\
\overline{t}&=& t+\epsilon\tau(x, t, u)+O(\epsilon^2),\\\label{eq:3}
\overline{u}&=& u+\epsilon\varphi(x, t, u)+O(\epsilon^2),\nonumber
\ee
where
$\displaystyle{\xi=\frac{{\p}\overline{x}}{{\p}\epsilon}\big|_{\epsilon=0},
~\tau=\frac{{\p}\overline{t}}{{\p}\epsilon}\big|_{\epsilon=0}}$,
  and
 $\displaystyle{\varphi=\frac{{\p}\overline{u}}{{\p}\epsilon}\big|_{\epsilon=0}}$.
The general vector field
\be
\X=\xi(x,t,u)\frac{\p}{{\p}x}+
\tau(x,t,u)\frac{\p}{{\p}t}+\varphi(x,t,u)\frac{\p}{{\p}u}~.\label{eq:infinitesimal}
\ee
on the $(x,t,u)$ space is assumed. We define the characteristic function
$Q=\varphi-\xi u_x-\tau u_t$. Then the third order
prolongation of the infinitesimal operator (\ref{eq:infinitesimal})
can be showed by the following prolongation formulas:
\be
{\rm Pr}^{(3)}\X=\xi(x,t,u)\frac{\p}{{\p}x}+
\tau(x,t,u)\frac{\p}{{\p}t}+\sum_{\#J=j=0}^3
\varphi_J(x,t,u^{(j)})\frac{\p}{{\p}u_J},\label{eq:6}
\ee
with coefficients
\be
\varphi_{J}=D_JQ+\xi u_{J,x}+\tau u_{J,t}.\label{eq:7}
\ee
Here, $J=(j_1, j_2, j_3)$
is a multi-indices, and $D_i$ is total derivative.

Using theorem \ref{thm:2} and relation (\ref{eq:inv}), we have
\be
{\rm Pr}^{(3)}\X\left[u_t+[a(t)u^n-l(t)]u_x+b(t)u_{xxx}-c(t)u\right]=0,
\ee
whenever
\bee
u_t+[a(t)u^n-l(t)]u_x+b(t)u_{xxx}-c(t)u=0.
\eee
Since $\xi, \tau$ and $\varphi$ only depend on $x, t, u$ one may
calculate the  coefficients to zer which  leads to the following
determining equations:
\be\label{eq:determinequation}
\left\{ \begin{array}{lcl}
\xi_u=\tau_u=\tau_x=0,\\[2mm]
\varphi_{xu}=\varphi_{uu}=0,\\[2mm]
3b(t)\xi_x=\dot{b}(t)\tau+b(t)\tau_t, \\[2mm]
\varphi_t+[(\varphi_u-\tau_t)c(t)-\tau\dot{c}(t)]u+\varphi_x a(t)u^n
-\varphi_x l(t)
-\varphi c(t)+\varphi_{xxx}b(t)=0,\\[2mm]
n\varphi a(t)u^{n-1}+\xi_xl(t)-\xi_uc(t)u-b(t)\xi_{xxx}-\tau\dot{l}(t)
-\xi a(t)u^n
+\tau_t(a(t)u^n-l(t))\\
-\xi_t+\tau\dot{a}(t)u^n=0.
\end{array}\right.
\ee
The general solution to system of partial differential equations (\ref{eq:determinequation}) is
\be\nonumber
\xi(x, t, u)&=&{1\over3}\left(\frac{\dot{b}(t)\beta(t)}{b(t)}+\dot{\beta}(t)\;\right)x+\alpha(t),\\
\tau(x, t, u)&=&\beta(t),\\\nonumber
\varphi(x, t, u)&=&\gamma(t)u+\eta(x, t),
\ee
where $b, \alpha, \beta, \gamma$ are arbitrary smooth functions with respect to $t$ and
$\eta$ also is a smooth function with respect to $x, t$.
\section{The Cylindrical KdV Equation}
One of a special case of (\ref{eq:maineq}) equation is $a(t) = 6, n = 1,
b(t) = 1$ with $c(t) = 0$ and $l(t)=-1/2t$  which
the (\ref{eq:maineq}) reduces to \cite{[2]}
\be\label{eq:cykdveq}
u_t+6uu_x+u_{xxx}+{1\over2t}u=0,
\ee
Using the (\ref{eq:infinitesimal}) vector field and its third prolong (\ref{eq:6}),
we have
\be
{\rm Pr}^{(3)}\X\left[u_t+6uu_x+u_{xxx}+{1\over2t}u\right]=0,\label{eq:cykv}
\ee
whenever
\bee
u_t+6uu_x+u_{xxx}+{1\over2t}u=0.
\eee
Solving (\ref{eq:cykv}) leads to following determining system
\be\label{eq:determinsys}
\left\{ \begin{array}{lcl}
\xi_u=\tau_u=\tau_x=0,\\[2mm]
\xi_{xx}=\varphi_{xu}=\varphi_{uu}=0,\\[2mm]
\tau_t=3\xi_x, \\[2mm]
[t\varphi_u-t\tau_t+12t^2\varphi_x-\tau]u+2t^2(\varphi_t+\varphi_{xxx})-t\varphi=0\\[2mm]
6\varphi-6u(\xi_x-\tau_t)u-\xi_t=0.
\end{array}\right.
\ee
By solving (\ref{eq:determinsys}) system with respect to $\xi, \tau$ and $\varphi$,
we obtain
\be\label{eq:desolKdV}
\xi(x, t)={c_1\over3}x+c_3t\sqrt{t}+c_2,
\qquad \tau(t)=c_1t,\qquad \varphi(x, t, u)=-{2\over3}c_1u+{c_3\over4}\sqrt{t}.
\ee
Therefore the infinitesimal generators are
\begin{subequations}\label{eq:infiniKdV}
\begin{align}
&\X_1={x\over3}{\p}_x+t{\p}_t-{2\over3}u{\p}_u,\\
&\X_2={\p}_x,\\
&\X_3=t\sqrt{t}{\p}_x+{\sqrt{t}\over4}{\p}_u,
 \end{align}
\end{subequations}
with following commutation relations
\be\label{eq:Liebracket}
[\X_1,\X_2]=-{1\over 3}\X_2,\quad [\X_1,\X_3]=-{7\over 6}\X_3,\quad [\X_2,\X_3]=0.
\ee
%
\section{Optimal System}
Assume $G$ is a Lie group and $\mathfrak{g}$ its Lie algebra.
For any element $T\in G$ we have a inner automorphism with definition
$T_a\longmapsto TT_aT^{-1}$ on the Lie group $G$. This
automorphism of the group $G$ induces an automorphism of $\mathfrak{g}$.
The group of all these automorphisms forms a Lie group that is called {\it
the adjoint group $G^A$}. For arbirary  $X,Y\in g$, we can define the linear mapping ${\rm
Ad}X(Y):Y\longrightarrow[X,Y]$ which is an automorphism of $\mathfrak{g}$,
called {\it the inner derivation of $\mathfrak{g}$}. For all $X,Y\in\mathfrak{g}$, the algebra of all
inner derivations ${\rm ad}X(Y)$ together with
the Lie bracket $[{\rm Ad}X,{\rm Ad}Y]={\rm Ad}[X,Y]$ is a Lie
algebra $\mathfrak{g}^A$ called the {\it adjoint algebra of $\mathfrak{g}$} which  $\mathfrak{g}^A$ is the Lie algebra of $G^A$.
Two subalgebras in $\mathfrak{g}$ are {\it conjugate}
if there is a transformation of $G^A$ which takes one subalgebra
into the other. The collection of pairwise non-conjugate
$s$-dimensional subalgebras is the optimal system of subalgebras
of order $s$. The construction of the one-dimensional optimal
system of subalgebras can be carried out by using a global matrix
of the adjoint transformations as suggested by Ovsiannikov
\cite{[4]}. The latter problem, tends to determine a list (that is
called an {\it optimal system}) of conjugacy inequivalent
subalgebras with the property that any other subalgebra is
equivalent to a unique member of the list under some element of
the adjoint representation i.e. $\overline{\mathfrak{h}}\,{\rm
Ad(g)}\,\mathfrak{h}$ for some ${\rm g}$ of a considered Lie group.
Thus we will deal with the construction of the optimal system of
subalgebras of $\mathfrak{g}$.
The adjoint action is given by the Lie series
\begin{eqnarray}
{\rm Ad}(\exp(s\,\X_i))\X_j=\X_j-s\,[\X_i,\X_j]+\frac{s^2}{2}\,[\X_i,[\X_i,\X_j]]-\cdots,
\end{eqnarray}
where $s$ is a parameter and $i,j=1,\cdots,n$.

We can expect to simplify a given arbitrary element,
\begin{eqnarray}\label{eq:vectorfield}
\X=\sum_{i=1}^3a_i\X_i,
\end{eqnarray}
of the Lie algebra $\mathfrak{g}$. Note that the elements of $\mathfrak{g}$ can be represented by vectors
$(a_1, a_2, a_3)\in{\Bbb R}^3$ since
each of them can be written in the form (\ref{eq:vectorfield}) for some constants $a_1, a_2, a_3$.
Hence, the adjoint action can be regarded as (in
fact is) a group of linear transformations of the vectors $(a_1, a_2, a_3)$.
%
\begin{thm}\label{thm:1}
An optimal system of one--dimensional Lie subalgebras of (\ref{eq:infiniKdV}) equation is
 generated by
\bee
(1)\quad A^1_1=\langle a\X_2+b\X_3\rangle  \qquad\qquad\qquad (2) \quad A^2_1=\langle \X_1\rangle
\eee
where $a, b\in{\Bbb R}$ are arbitrary constants.
\end{thm}
\begin{proof}
Suppose that
$F^{\varepsilon}_i:{\mathfrak g}\rightarrow{\mathfrak g}$ defined by
$\X\mapsto\mbox{Ad}(\exp(\varepsilon \X_i)\X)$
is a linear map, for $i=1,2,3$. The matrices
$M^\varepsilon_i$ of $F^\varepsilon_i$, $i=1,2,3$, with
respect to basis $\{\X_1,\X_2,\X_3\}$ are
\bee
\begin{array}{lcl}
M^\varepsilon_1=\left[\begin{array}{ccccccc}
1&0&0\\
0&\exp({1\over3}\varepsilon)&0\\
0&0&\exp(-{7\over6}\varepsilon)
\end{array}
\right],\quad
M^\varepsilon_2=\left[\begin{array}{ccccccc}
1&-{1\over3}\varepsilon&0\\
0&1&0\\
0&0&1
\end{array}
\right],\quad
M^\varepsilon_3=\left[\begin{array}{ccccccc}
1&0&{7\over6}\varepsilon\\
0&1&0\\
0&0&1
\end{array}
\right].
\end{array}
\eee
Let $\displaystyle{\X=\sum_{i=1}^3a_i\X_i}$, then we have
\bee
F_3^{s_3}\circ F_2^{s_2}\circ F_1^{s_1}: \X\mapsto[a_1-{s_2\over3}a_2+{7s_3\over6}a_3]\X_1 +[\exp({s_1\over3})a_2]\X_2+[\exp(-{7s_1\over6})a_3]\X_3.
\eee
If $a_2,a_3\neq0$ then we can omit the coefficients of $\X_1$  by setting
$\displaystyle{s_2 = {3a_1\over a_2}}$ and $\displaystyle{s_3 = -{6a_1\over 7a_3}}$
 So, $\X$ is reduced to the case (1). But if $a_2,a_3=0$, then $\X$ is reduced to the case (2).
 There is no any new case.
\end{proof}
%
\section{Symmetry Reductions and Exact Solutions}
The invariants associated with the infinitesimal generator $\X_2$ are obtained by integrating the characteristic equation
\be
{3dx\over x}={dt\over t}={-3du\over 2u},
\ee
which generates the invariants
\be\label{solve1}
r={x^3\over t},\quad g(r)=x^2u(x, t).
\ee
Substituting (\ref{eq:inv2}) into (\ref{eq:cykdveq}), to determine the form of the function $g$, the (\ref{eq:cykdveq})
is reduced to following third order ordinary differential equation
\be\label{invsol:1}
54r^3\dddot{g}-2r^2\dot{g}+84rg(r)\dot{g}-24g^2(r)-(48+r)g(r)=0,
\ee
with respect to $g(r)$ and here $\displaystyle{\dot{g}={dg\over dr}}$ and  $\displaystyle{\dddot{g}={d^3g\over dr^3}}$.

The characteristic equation associated with $\X_3$ is
\be
{dt\over t\sqrt{t}}={4du\over \sqrt{t}},
\ee
which generates the invariants $x, \sqrt[4]{t}\exp(u)$. Then the similarity solution have the form
\be
u(x, t)=\ln\left({f(x)\over\sqrt[4]{t}}\right).\label{eq:inv2}
\ee
By substituting (\ref{eq:inv2}) into (\ref{eq:cykdveq}), to determine the form of the function $f$, the (\ref{eq:cykdveq})
is reduced to following third order ordinary differential equation
\be\label{invsol:2}
4tf^2\dddot{f}-12tf\dot{f}\ddot{f}+8t\dot{f}^3+24t\ln({f\over\sqrt[4]{t}})f^2\dot{f}
-2\ln({f\over\sqrt[4]{t}})f^3-f^3=0,
\ee
with respect to $f$ and here $\displaystyle{\dot{f}={df\over dx}}$.

The associated  invariants of  $\X_2$ are the
arbitrary function $h(t, u)$.
%
\section{Numerical Analysis}
In this section, we use Chebyshev pseudo-spectral method (CPSM) to solve the
introduced problems (\ref{invsol:1}) and (\ref{invsol:2}). This method applies the Chebyshev-Gauss-Lobatto points
$$\xi _{j} =\cos\left(\frac{j\pi }{N} \right)\, ,\, \, j=0,\ldots ,\, N$$
as collocation points, that satisfy $T'(\zeta _{j} )(1-\zeta _{j} ^{2} )=0$ where
$T_N(x)$ is the Chebyshev polynomial of degree  $N$. Then, the Lagrange interpolating polynomials based on
$\xi _{j} \, ,\, \, j=0,\ldots ,\, N$ can be got as follows:
\be
L_{N,j} (x)=\frac{(-1)^{j+1} (1-x^{2} )T'_{N} (x)}{c_{j} N^{2} (x-\zeta _{j} )} \qquad j=0,\ldots ,\, N,
\ee
where
$\displaystyle{c_{j} =\left\{\begin{array}{l} {2\, \, \, \, j=0,N} \\ {1\, \, \, \, \, j=1,\ldots ,N-1} \end{array}\right.}$.
It is clear that $L_{N,j} (\zeta _{k} )=\delta _{j\, k}$ where $\delta _{j\, k}$ denotes Kronecker delta.
Therefore, a function  $z(x)$ is approximated in interval  $[-1, 1]$ as below:
\be
z(x)\approx \sum _{j=0}^{N}L_{N,j} (x)\, z(\zeta _{j} ).
\ee
Also, we can obtain  approximation for derivative values at collocation points
$\zeta _{i} \, (i=0,\ldots ,N)$ for  $z(x)$ as follows:
\be
z'(\zeta _{i})\approx \sum _{j=0}^{N}L'_{N,j} (\zeta _{i})z(\zeta _{j}) =\sum _{j=0}^{N}d_{i\, j} z(\zeta _{j}) \qquad  i=0,\ldots ,N,
\ee
where
$\displaystyle{D_N=\left(d_{ij}\right)_{i,j=0}^N}$
denotes Chebyshev collocation derivative matrix with
$\displaystyle{d_{i,j}=L'_{N,j}(\zeta_i)}$
 that can be got as \cite{[Peyret]}:
 \be
 c_{j} =\left\{
 \begin{array}{l}
 {\displaystyle{d_{ij}=\frac{c_i(-1)^{i+j}}{c_j(\zeta_i-\zeta_j)}} \qquad i,j=0,\cdots, N\;{\rm and}\; i\neq j} \\ [2mm]
 {\displaystyle{d_{ii}=\frac{-\zeta_i}{2(1-\zeta_i^2)}} \qquad\; i=1, \cdots, N-1,}\\[2mm]
 {\displaystyle{d_{0,0}=-d_{N,N}=\frac{2N^2+1}{6}}}
 \end{array}\right.
 \ee
So, in the matrix form, we can write
${\bf z}_1=D_N{\bf z}$
where
${\bf z}=\left[z(\zeta_0), \cdots, z(\zeta_N)\right]^T$
and
${\bf z}'=\left[z'(\zeta_0), \cdots, z'(\zeta_N)\right]^T$.
Now, in a similar way,
$z^{(k)}(\zeta_j)$ for $j=0, \cdots, N$
can be approximated by
${\bf z}_k\approx D_N^k {\bf z}$
with
${\bf z}_k=\left[z^{(k)}(\zeta_0), \cdots, z^{(k)}(\zeta_N)\right]^T$
where
$\displaystyle{D_N=\left(d_{ij}^{(k)}\right)_{i,j=0}^N}$
represents the
$k$-th
 power of $D_N$. Note that  $d_{ij}=d_{ij}^{(1)}$. So we have
 \be\label{eq:iii}
 z^{(k)}(\zeta_i)\approx \sum_{j=0}^N d_{ij}^{(k)} z(\zeta_j), \qquad i=0, \cdots, N.
 \ee
\subsection{Problem 1: CPSM for (\ref{invsol:1})}
\begin{figure}
\includegraphics[height=5cm]{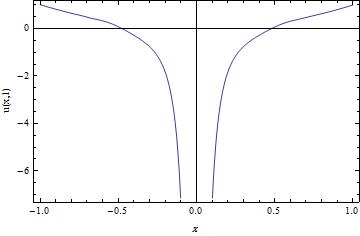}
\includegraphics[height=5cm]{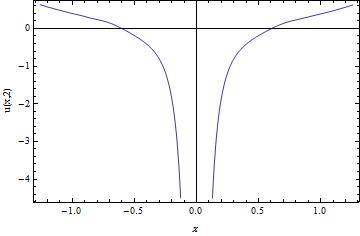}
\centering{\includegraphics[height=5cm]{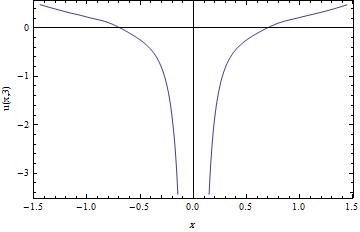}}
\caption{Plot of approximate solution $u(x, t)$  for $t=1,2,3$ and $N=25$ in Problem 1.} \label{fig1}
\end{figure}

We apply CPSM for solving the differential equation (\ref{invsol:1}) with boundary conditions $g(-1)=g'(-1)=g(1)=1$.
 By employing the approximate formulas of derivatives (\ref{eq:iii}), the problem is reduced as below:
 \bee
 \left\{
 \begin{array}{l}
 {\displaystyle{54\zeta_i^3\left(\sum_{j=0}^Nd_{ij}^{(3)}g(\zeta_j)\right)+(84\zeta_ig(\zeta_i)-2\zeta_i^2)\left(\sum_{j=0}^Nd_{ij}g(\zeta_j)\right)+24g^2(\zeta_i)}}\\
    {\displaystyle{-(48+\zeta_i)g(\zeta_i)=0,} \qquad i=0, \cdots, N-2,} \\ [2mm]
 {\displaystyle{g(\zeta_N)=g(\zeta_0)=1, \quad \sum_{j=0}^Nd_{Nj}g(\zeta_j)=1}}
 \end{array}\right.
 \eee
Therefore, we have a nonlinear system of
 $N-1$
 equations and  $N-1$ unknown parameters  $g(\zeta_i)$ , for $i=1, \cdots, N-1$,
 which can be solved by Newton's method.
 We set the obtained approximate solutions for collocations points $\zeta_i$, for $i=1, \cdots, N-1$,
  in the problem (\ref{invsol:1}) and get the residuals for these points (i.e.  $L[g(\zeta_i)]$,
  if $L$ be the operator of the problem (\ref{invsol:1}) that it operate on function  $g$).
   We can see the results for $N=25$  in Table \ref{table:nonlin1}.
   Also, from (\ref{solve1}), we can observe the behaviors of solutions $u(x, t)$  for  $t=1, 2, 3$ in Figures \ref{fig1}.
   The results show that the obtained solutions have high accuracy.
\begin{table}[ht]
\caption{Residuals for $N=25$ in Problem 1.}
\centering
\begin{tabular}{c c}
\hline\hline
  $i$ & $|L[g(\zeta_i)]|$ \\[2mm]
  \hline
  1 & $4.97523869142924*10^{-5}$\\[2mm]
2& $7.23248462008996*10^{-6}$\\[2mm]
3& $8.63460245170699*10^{-6}$\\[2mm]
4& $1.49965897122683*10^{-6}$\\[2mm]
5& $1.17603000404642*10^{-6}$\\[2mm]
6& $6.85683971823891*10^{-8}$\\[2mm]
7& $7.78958764158233*10^{-9}$\\[2mm]
8& $3.41905703749034*10^{-9}$\\[2mm]
9& $5.87228043968934*10^{-11}$\\[2mm]
10& $1.35065292283798*10^{-11}$\\[2mm]
11& $3.37907479774912*10^{-12}$\\[2mm]
12& $1.07969189144796*10^{-13}$\\[2mm]
13& $1.16351372980716*10^{-13}$\\[2mm]
14& $2.84927637039800*10^{-12}$\\[2mm]
15& $1.23083765402043*10^{-11}$\\[2mm]
16& $1.10974340827851*10^{-10}$\\[2mm]
17& $2.13702122664471*10^{-9}$\\[2mm]
18& $1.56747148594149*10^{-9}$\\[2mm]
19& $1.25001037076799*10^{-7}$\\[2mm]
20& $1.83994266933495*10^{-7}$\\[2mm]
21& $1.65322587974969*10^{-6}$\\[2mm]
22& $4.29589790940099*10^{-6}$\\[2mm]
23& $2.31182675776153*10^{-5}$\\[2mm]
24& $1.02044087658533*10^{-6}$\\[2mm]
  \hline
\end{tabular}\label{table:nonlin1}
\end{table}
\begin{table}[ht]
\caption{Residuals for $N=25$ in Problem 2.}
\centering
\begin{tabular}{c c c c c c}
\hline\hline
  $i$ & $t=1$& $t=2$& $t=3$\\[2mm]
  \hline
1 & $1.079968190431*10^{-7}$ &  $2.16307934230997*10^{-7}$ & $3.328058753027108*^{-7}$  \\[2mm]
2 & $1.355961942728*10^{-8}$ &  $2.81342655839011*10^{-8}$ & $4.211412374388601*^{-8}$  \\[2mm]
3 & $3.660495906387*10^{-9}$ &  $7.69985497584002*10^{-9}$ &   $1.167282093206267*^{-8}$\\[2mm]
4 & $1.316350584090*10^{-9}$ & $2.77901524015078*10^{-9}$ &  $4.188224878021174*^{-9}$ \\[2mm]
5 & $4.581313106655*10^{-10}$ & $9.81683345724349*10^{-10}$ &  $1.433773988424036*^{-9}$ \\[2mm]
6& $1.671205396291*10^{-10}$ & $2.87156076694827*10^{-10}$ &   $4.428528654898400*^{-10}$\\[2mm]
7& $6.472511415722*10^{-11}$ & $1.17186260695234*10^{-10}$ &   $1.448259290270925*^{-10}$ \\[2mm]
8& $6.582823175449*10^{-11}$ & $1.23897336834488*10^{-10}$ &   $1.694715479061415*^{-10}$ \\[2mm]
9& $8.133582696245*10^{-11}$ & $1.53065116137440*10^{-10}$ &   $1.938946780910555*^{-10}$ \\[2mm]
10& $4.537525910564*10^{-11}$ & $1.10620845816811*10^{-10}$ &   $1.189839338167075*^{-10}$ \\[2mm]
11& $4.972022793481*10^{-11}$ & $6.80735467994964*10^{-11}$ &   $9.872991313386592*^{-11}$ \\[2mm]
12& $1.749356215441*10^{-11}$ & $4.08562073062057*10^{-13}$ &   $2.186339997933828*^{-11}$ \\[2mm]
13& $2.106403940160*10^{-11}$ & $3.94928534319660*10^{-11}$ &   $4.613376347606390*^{-11}$ \\[2mm]
14& $2.006395050102*10^{-12}$ & $1.84368076361352*10^{-11}$ &   $1.365929591656822*^{-11}$ \\[2mm]
15& $3.616218435809*10^{-11}$ & $5.69819746942812*10^{-11}$ &   $5.420552895429864*^{-11}$ \\[2mm]
16& $4.798295094587*10^{-11}$ & $1.01785246897634*10^{-10}$ &   $1.384812264859647*^{-10}$ \\[2mm]
17& $6.135358887604*10^{-11}$ & $1.11201714503295*10^{-10}$ &   $1.379252267952324*^{-10}$ \\[2mm]
18& $5.050537765782*10^{-11}$ & $1.04055430938387*10^{-10}$ &   $1.556180739825663*^{-10}$ \\[2mm]
19& $2.251852038170*10^{-10}$ & $5.15020026625734*10^{-10}$ &   $7.070202201475695*^{-10}$ \\[2mm]
20& $5.082538834244*10^{-10}$ & $1.06219899542736*10^{-9}$ &   $1.523963177874065*^{-9}$ \\[2mm]
21& $1.287500328572*10^{-9}$ & $2.73905342851321*10^{-9}$ &    $3.971750928144502*^{-9}$ \\[2mm]
22& $3.731360831427*10^{-9}$ & $7.70869235111604*10^{-9}$ &   $1.125179949212906*^{-8}$ \\[2mm]
23& $1.381041303538*10^{-8}$ & $2.82384888805609*10^{-8}$ &   $4.212327198160892*^{-8}$ \\[2mm]
24& $2.073358654819*10^{-5}$ & $6.92459967588376*10^{-6}$ &    $3.196776287239799*^{-6}$ \\[2mm]
  \hline
\end{tabular}\label{table:nonlin2}
\end{table}
\begin{figure}
\includegraphics[height=5cm]{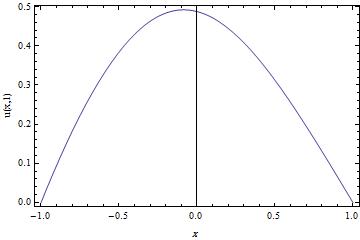}
\includegraphics[height=5cm]{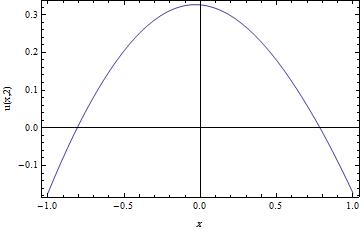}
\centering{\includegraphics[height=5cm]{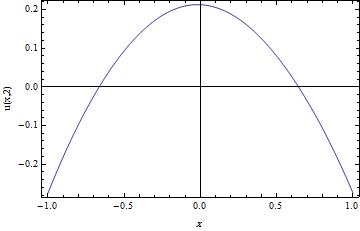}}
\caption{Plot of approximate solution $u(x, t)$  for $t=1,2,3$ and $N=25$ in Problem 2.} \label{fig2}
\end{figure}
%
\subsection{Problem 2: CPSM for (\ref{invsol:2})}
We can use CPSM for the differential equation (\ref{invsol:2}) that obtained the previous section.
We consider boundary conditions  $f(-1)=f'(-1)=f(1)=1$ for this problem.
After inserting the formulas (\ref{eq:iii}), we have the reduce form:
 \bee
 \left\{
 \begin{array}{l}
 {\displaystyle{4tf^2(\zeta_i)\left(\sum_{j=0}^Nd_{ij}^{(3)}f(\zeta_j)\right)
 -12tf(\zeta_i)\left(\sum_{j=0}^Nd_{ij}f(\zeta_j)\right)\left(\sum_{j=0}^Nd_{ij}^{(2)}f(\zeta_j)\right)
  +8t\left(\sum_{j=0}^Nd_{ij}f(\zeta_j)\right)^3}}\\
  {\displaystyle{+24tf^2(\zeta_i)\ln\left(\frac{f(\zeta_i)}{\sqrt[4]{t}}\right)\left(\sum_{j=0}^Nd_{ij}f(\zeta_j)\right)
  -2f^3(\zeta_i)\ln\left(\frac{f(\zeta_i)}{\sqrt[4]{t}}\right)-f^3(\zeta_i)=0,}\qquad i=0, \cdots, N-2.}\\[5mm]
  {\displaystyle{f(\zeta_N)=f(\zeta_0)=1, \quad \sum_{j=0}^Nd_{Nj}f(\zeta_j)=1.}}
 \end{array}\right.
 \eee
Similar to the Problem 1, for any $t$, we get a system of $N-1$ nonlinear equations and  $N-1$ unknown parameters $f(\zeta_i)$ for $i=1, \cdots, N-1$ .
Accuracy of the approximate solutions for $N=25$ and $t=1,2,3$  is observable in Table \ref{table:nonlin2}.
Also, by (\ref{eq:inv2}) we can show the behaviors of $u(x, t)$  for $N=25$ and  $t=1, 2, 3$  in Figures \ref{fig2}.

%
\section{Conclusion and Results}
Lie point symmetries of the Korteweg-de Vries equation with power-law nonlinearity (\ref{eq:maineq})
in a particular form of  is $a(t) = 6, n = 1,
b(t) = 1$ with $c(t) = 0$ and $l(t)=-1/2t$,
form a three dimensional Lie symmetry algebra.
The invariant solution of these symmetries is reduced
The optimal system of one-dimensional Lie subalgebras associated to this symmetry algebra
 is generated by two vector fields.

 The Chebyshev pseudo-spectral method (CPSM) as a numerical analysis
  is applied for invariant solutions. In this method, we use the Chebyshev-Gauss-Lobatto points
as collocation points.

\end{document}